% Manuscript of the paper Analyticity of functions analytic on circles
% Written by Josip Globevnik
% Submitted for publication to Journal of Mathematica Analysis and Application on May 7, 2009
% Accepted on June 4, 2009
\magnification 1200
\def\R{{\rm I\kern-0.2em R\kern0.2em \kern-0.2em}}
\def\N{{\rm I\kern-0.2em N\kern0.2em \kern-0.2em}}
\def\P{{\rm I\kern-0.2em P\kern0.2em \kern-0.2em}}
\def\B{{\rm I\kern-0.2em B\kern0.2em \kern-0.2em}}
\def\Z{{\rm I\kern-0.2em Z\kern0.2em \kern-0.2em}}
\def\C{{\bf \rm C}\kern-.4em {\vrule height1.4ex width.08em depth-.04ex}\;}

\def\D{{\Delta}}
\def\DD{{\overline \Delta}}
\def\bD{{b\Delta}}
\def\z{{\zeta}}

\def\cS{{\cal S}}
\def\cS{{\cal S}}

\def\G{{\Gamma}}

\font\ninerm=cmr8
\
\vskip 22mm
\centerline {\bf ANALYTICITY OF FUNCTIONS ANALYTIC ON CIRCLES}
\vskip 8mm
\centerline{Josip Globevnik}
\vskip 8mm
{\noindent \ninerm ABSTRACT \ \ Let $\D$ be the open unit disc in $\C$, let $p\in \bD, $ and
let $f$ be
a continuous function on $\DD$ which extends holomorphically from each circle in $\DD$
centered at the origin and from each circle in $\DD$ which passes through $p$. Then
$f$ is holomorphic on $\D$.}
\vskip 6mm
\bf 1.\ Introduction and the main result
\vskip 2mm
\rm Denote by $\D$ the open unit disc in $\C$. If $f$ is a continuous function
on a circle $\G $ then we say that $f$ extends
holomorphically from $\Gamma $ if it extends holomorphically through the disc
bounded by $\G$.

Let $f$ be a continuous function on $\DD$ which extends holomorphically from
every circle $|\z|=r,\ 0<r\leq 1$. A trivial example is a function constant on each
circle $|\z | =r$. Obviously
such a function is not necessarily holomorphic on $\D$. There are worse
examples. For instance, the function
$$
f(z) = \left\{ \eqalign{ &z^2/\overline z\ \ \ (z\in \DD\setminus \{ 0\} \cr &0\ \ \
\ \ \ \ (z=0)\cr}\right.
$$ is  continuous on $\DD $ and extends holomorphically from every circle $\G$ in
$\DD$ that either
surrounds the origin or contains the origin, yet $f$ is not holomorphic on $\D$ \ [G2].

Let $p\in\bD$. In the present paper we show that if a continuous function on $\DD$
extends holomorphically
from each circle centered
at the origin and also from each circle contained in $\DD$ and passing through $p$ then
it must be holomorphic on $\D$. In fact, we prove a somewhat better result:
\vskip 2mm
\noindent\bf THEOREM 1.1\ \it Let $p\in\bD$ and let $\tau < 1/2$. Suppose that $f$ is
a continuous function on $\DD$ such that

(i)\ $f$ extends holomorphically from each circle $|\z |=R, \ 0<R\leq 1$

(ii)\ $f$ extends holomorphically from each circle of radius $R\geq \tau $ which is
contained in
$\DD$ and passes through $p$.

\noindent Then $f$ is holomorphic on $\D$. \rm
\vskip 2mm
For each $z\in\C,\ r>0$, denote $\D (z,r)=\{\z\in\C\colon\ |\z - z|<r\}$. Our family
of circles can be written as $\{
b\D (a(t),r(t))\colon \ 0<t<1\} $ where $t\mapsto a(t),\ t\mapsto r(t)$ are piecewise
smooth functions on $[0,1]$. Tumanov [T2] proved that continuous functions that extend
holomorphically from each
circle belonging to such a family are holomorphic provided that $\DD (a(0),r(0))\cap
\DD (a(1),
r(1))=\emptyset$ and provided that no circle $\bD (a(s), r(s))$ is contained in the
closed disc $\DD (a(t),r(t))$ if $t\not = s$. Note that
the second condition is not satisfied by our family. More general results are known
in real-analytic category [A]. Note that our family of circles is not real-analytic.
\vskip 4mm
\bf 2.\ Semiquadrics and the related problem in $\C ^2$\rm
\vskip 2mm
We begin the proof of our theorem. With no loss of generality assume that $p=-1$. As
in [AG] and [G1] we introduce semiquadrics
to pass to an associated problem in $\C^2$. Given $a\in\C$ and $r>0$ let
$$
\Lambda _{a,r} = \{ (z,w)\in\C^2\colon\ (z-a)(w-\overline a)=r^2,\ 0<|z-a|<r\} .
$$
This is a closed complex submanifold of $\C^2\setminus\Sigma $ where
$\Sigma = \{ (\z, \overline\z )\colon\ \z\in\C\}$, which is attached to $\Sigma $ along
$b\Lambda _{a,r}= \{ (\z,\overline \z)\colon\ \z\in\bD (a,r)\}$. A
continuous function $g$ extends holomorphically from the circle $\bD (a,r)$ if and
only if the
function $G$, defined on $b\Lambda_{a,r}$ by
$G(\z,\overline\z )= f(\z )\ (\z\in b\D (a,r))$ has a bounded continuous extension
to $\overline {\Lambda _{a,r}} =
\Lambda _{a,r}\cup b\Lambda _{a,r}$ which is holomorphic on $\Lambda _{a,r}$. In fact,
if we denote by the same letter $g$
the holomorphic extension of $g$ through $\D (a,r)$ we have
$$
G\Bigl(z,\overline a + {{r^2}\over{z-a}}\Bigr) = g(z)\ \ (z\in\DD(a,r)\setminus \{ a\})
$$
and, if we define $G(a,\infty)=g(a)$ we get a continuous function $G$ on
$\tilde\Lambda _{a,r}
= \overline {\Lambda _{a,r}}\cup \{(a,\infty)\}$, the closure of $\Lambda_{a,r}$ in
$\C\times \overline\C$.

It is known that if $(a,r)\not= (b,\rho)$ then $\Lambda _{a,r}$
meets $\Lambda_{b,\rho}$ if and only if $a\not= b$ and one of the circles
$b\D (a,r),\ b\D(b,\rho)$ surrounds the other [G1].

Let $\tau$ and $f$ be as in Theorem 1.1. By our assumption, $f$ extends
holomorphically from two families of circles: $\{ b\D (t,t+1):\ -1+\tau\leq t\leq 0\}$ and
$\{ b\D (0,R):\ 0<R\leq 1\}$. Accordingly, there are two families of semiquadrics:\
$\{ \Lambda _{t,t+1}:\ -1+\tau \leq t\leq 0\} $ and  $\{ \Lambda_{0,R}:\ 0<R\leq 1\}$ and the function
$F(\z,\overline \z) = f(\z)\ (\z\in\D)$ has a bounded holomorphic
extension through each of these semiquadrics.
In each of these families the semiquadrics are pairwise disjoint. Let us look first at
the first family and let $N$ be the closure of the union of
$\Lambda _{t,t+1},\ -1+\tau\leq t\leq 0$ in $\C\times\overline\C$, that is,
$$
N = \bigcup _{-1+\tau\leq t\leq 0} \bigl[ \Lambda _{t,t+1}\cup b\Lambda _{t,t+1}\cup
\{ (t,\infty)\} \bigr] .
$$
The contuinuity of $f$ together with the maximum principle implies that
our functon $(\z ,\overline\z)\mapsto F(\z,\overline\z ) = f(\z )$ defined on $\{
(\z,\overline\z)\colon\ \z\in\D\} $
extends from $N\cap\Sigma =
\{ (\z,\overline\z)\colon\ \z\in \DD\setminus \D(-1+\tau ,\tau)\} $ continuously
to $N$ so that the extension $F$
is holomorphic on each fiber $\Lambda _{t,t+1},\ -1+\tau\leq t\leq 0$. Note that
the part $N_0$ of $N$ contained in $\C\times\C$ is a smooth
CR manifold with piecewise smooth boundary consisting of three
smooth pieces:\ $\overline{\Lambda _{-1+\tau,\tau}},\ \overline{\Lambda_{0,1}}$ and
$N\cap\Sigma $ and the function $F$ is CR in the interior, that is,
$$
\int _{N_0} f\overline\partial\omega = 0
$$
for each smooth $(2,0)$-form $\omega $ on $\C^2$ whose support intersects the interior of $N_0$ in
a compact set.

Now look at the second family and let $L$ be the closure
of the union of $\Lambda _{0,R},\ 0<R\leq 1$, in $\C\times \C$, that is
$$
L= [\{ 0\}\times\C]\cup [\cup_{0<R\leq 1}\overline {\Lambda _{0,R}}] .
$$
Again, our function $F$ extends from $L\cap\Sigma = \{ (\z, \overline\z )\colon\
\z\in\DD\}$
to a bounded continous
function on $L$ which is holomorphic on each leaf $\Lambda _{R,0},\ 0<R\leq 1$.
Again, the part $L_0$ of $L$ contained in $\C^2\setminus (\{ 0\}\times\C)$
is a CR manifold with piecewise smooth boundary consisting of
two pieces:\ $\overline{\Lambda_{0,1}}$
and $\{ (\z,\overline\z)\colon\ \z\in \DD\setminus\{ 0\}\}$ and
the extension $F$ is CR on the interior of $L_0$.

Tumanov's condition that no circle $b\D (a(s),r(s))$ is contained in the
closed disc $\overline{\D (a(t),r(t))}$ if $s\not= t$ implies that the semiquadrics
$\Lambda _{a(t),r(t)}$
are pairwise disjoint so their union is a CR manifold through which the function $F$ extends
as a CR function. Tumanov then uses an argument of H.Lewy [L] and the
Liouville theorem to show that the function $F$ does not depend on the second
variable,
that is, that $f$ is
holomorphic. We want to follow the same idea but in our case the semiquadrics are no
more
pairwise disjoint and so their union is not a manifold.
 In particular,
the manifolds $L$ and $N$ intersect. However, we show that our particular
geometric setting allows to apply the reasoning of Tumanov
on  $(L\cup N)\setminus (L\cap N)$, a CR manifold to which $F$ extends
as a CR function, to be able to conclude that the function 
$F$ does not depend on the second variable. We provide a detailed proof of Theorem 1.1. 
\vskip 4mm
\bf 3.\ The manifolds $L$ and $N$
\vskip 2mm
\rm As we have already mentioned, the function $F$ extends to $L$ and to $N$ so
that the extensions are
holomorphic on semiquadrics, the holomorphic fibers of $L$ and $N$. There is one
piece of $L\cap N$,
namely $\Lambda_{0,1} $ on which both extensions coincide.
However, a semiquadric of $L$ can intersect a semiquadric of $N$. In
fact, $\Lambda _{0,R}$ intersects $\Lambda _{t,t+1}$ if and only if $R<2t+1$.
We know that in this case the intersection consists of one point [G1]. It is easy to
see that it
is of the form $(x,y)$ where $x>0$ and $y>0$. This implies that there
is no problem in defining the extension of $F$ in the part of $N\cup L$ which is
contained in
$(\C\times\overline\C)\setminus ([0,1]\times\overline\C)$. We denote this part of $L\cup N$ by
$M$.
Let $\cS = \D\setminus \bigl\{\overline
{\D (\tau,1+\tau)}\cup [0,1]\bigr\} $. Given $z\in\cS$ we shall study
$M_z = \{\z\in\C\colon\ (z, \z)\in M\} $. We shall show that if $\Im z\not= 0,\
M_z$ is a closed
curve consisting of the segment joining $\overline z$ and $1/z$ and a circular arc
joining $1/z$ and $\overline z$ and
if $z\in\R$ then $M_z$ is the real axis in $\C$.
\vskip 2mm
\noindent\bf LEMMA 3.1\ \it\ \ Let $z\in\cS,\ z\not\in \R$. The circle $C_z$ passing
through
$1/z,\ \overline z$ and ${-1}$ is tangent to the real axis at $-1$. Let $\lambda _z$ be the
arc of $C_z$ with end points $\overline z$ and $1/z$ which does not contain $-1$.
Then $M_z$ consists of $\lambda _z$ and of the segment joining $\overline z$ and $1/z$. \rm
\vskip 2mm
\noindent \bf Proof.\ \rm We have $L\cap (\{ z\}\times \C) = \{ (z,R^2/z)\colon\  |z|\leq R\leq 1\}$
which is the segment joining $(z,\overline z)$ and $(z,1/z)$. To find what $N\cap (\{ z\}\times \C)$
is we recall first that
$$
\Lambda _{t,t+1} = \Bigl\{ (z,w)\colon\ w=t+{{(t+1)^2}\over{z-t}},\ |z-t|<1+t\Bigr\}.
$$
So we must determine $\{ w(t)\colon\ t(z)\leq t\leq 0\}$ where
$$
w(t)= t+{{(t+1)^2}\over{z-t}} =  {{(z+2)t+1}\over{z-t}}
 $$
and where $t(z)$ is such that $z$ lies on the circle $|\z -t(z)|=t(z)+1$, that is, when
 $w(t(z))=\overline z$. To find what circle
$$
\Bigl\{ {{(z+2)t+1}\over{z-t}}\colon\ t\in\R \Bigr\}
\eqno (3.1)
$$
is, write $z=P+iQ $ with $P, Q$ real,  and assume that $Q\not= 0$. We have
$$\eqalign
{
   {{(P+iQ+2)t+1}\over{P+iQ-t}} &=˛\cr
   &={{ [(P+iQ)t+2t+1].[P-iQ-t]}\over{(P-t)^2+ Q^2}} \cr
   &= {{(P^2+Q^2)t+2t(P-iQ)+P-iQ- (P+iQ)t^2-2t^2-t}\over {(P-t)^2+ Q^2}} .\cr}
$$
This is real when
$$
0 = -2tQ-Q-Qt^2= Q(t+1)^2, 
$$
that is, when $t= -1$ when
$$
w={{(z+2)(-1)+1}\over {z-(-1)}} = -1.
$$
It follows that (3.1) is a circle tangent to the real axis at $-1$ and it also follows that the arc
from $\overline z$ to $1/z$ containing $-1$ does not belong to $\{ \z\colon\ (z,\z )\in N\} $.
We already know that this circle must contain $\overline z$ and $1/z$. This completes the proof.
\vskip 2mm
\noindent We may compute the center of the circle in Lemma 3.1 by intersecting the real
line $\{ (-1+\overline z)/2+i\lambda (-1-\overline z)/2\colon\ \lambda\in\R\}$
with the vertical line $\Re \z = -1$. Again, write $z=P+iQ$ with $P, Q$ real. We
compute $\lambda $ at the point of intersection from the condition  $(1/2)(-1+P-\lambda Q) = -1$
which gives $ \lambda = (P+1)/Q$ and a short computation
shows that the center of the circle is $-1-i|z+1|^2/(2\Im z)$.

We now look at what $M_T$ is when $T\in\cS $ is real, that is, when $-1+2\tau <T<0$. Observe first
that $\Lambda _{0,R}$ intersects $\{ T\}\times \C$ if $|T|<R<1$ and
$(\{T\}\times\C)\cap\Lambda_{0,R}= \{(T,R^2/T)\}$. When $R$ moves from $|T|$ to $1$ the point
$R^2/T$ moves on the real axis from $T$ to $1/T$.
 This takes care of the intersection of $\{ T\}\times \C$ with $L$.
To find the intersection with $N$ we have to see what
$$
w=t+{{(1+t)^2}\over{T-t}} ={{t(T+2)+1}\over{T-t}} = w(t)
$$
does when $t$ decreases from $0$ to $(T-1)/2$, that is, when $\Lambda_{t,t+1}$
meets $\{ T\}\times \C$. At $t=0$ we have $w(0)=1/T$ and as $t$ decreases from $0$
to $T$, $w(t)$ moves from
$1/T$ to $-\infty$ along the real axis. When $t$ decreases from $T$ to $(T-1)/2$, $w(t)$
decreases from $+\infty$ to $T$. Thus, $M_T=\R\cup\{\infty\}$.
\vskip 4mm
\bf 4.\ Completion of the proof of Theorem 1.1 \rm
\vskip 2mm
Denote by $\pi (z,w) = z$ the projection onto the first coordinate axis.
Let $U\subset \cS$ be a small open
disc and consider $\pi^{-1}(U)\cap M$. This set consists of two smooth manifolds
$\pi^{-1}(U)\cap N$ and $\pi ^{-1}(U)\cap L$ with common boundary
$\{(\z, 1/\z)\colon\ \z\in U\}\cup\{(\z,\overline\z)\colon\ \z\in U\}$ along which
they meet transversely.
The set $\pi^{-1}(U)\cap M$ is a topological manifold which
can be oriented as part of the boundary of $\cup _{z\in U}\{ z\}\times D_z$ where,
for each $z\in \cS$, $D_z$ is the domain
bounded by $M_z$. The domains $D_z$ change continuously with $z\in\cS\setminus\R$
and as $z$ approaches $a\in\bD\setminus\R$ they
shrink to the point $\overline a$.

Since  we want to provide a detailed proof of Theorem 1.1. we shall need
\vskip 2mm
\noindent\bf LEMMA 4.1\  \it Let $B$ be an open ball in $\C^2$ and let $E\subset B$ be a
closed two-dimensional smooth submanifold of $B$
which is the common boundary of two closed three dimensional smooth submanifolds $\Sigma_1$
and $\Sigma_2$ of $B\setminus E$ such that $M=\Sigma_1\cup\Sigma_2\cup E$ is a
topological submanifold of $B$. Let $f$ be a continuous function on $M$  which is CR
on $\Sigma_1$ and $\Sigma_2$, that is, $
\int _{\Sigma_i} f\overline\partial \alpha = 0$ for each smooth,
$(2,0)$ form on $B$ whose support intersects $\Sigma_i$ in
a compact set, $i=1,2$. Then $f$ is CR on $M$, that is, $
\int _M f\overline\partial \alpha = 0$ for every smooth, $(2,0)$ form on $B$
whose support intersects $M$ in a compact set. \rm
\vskip 2mm
\noindent\bf Proof.\ \rm The proof, suggested by E.L.Stout, uses the fact obtained
by G.\ Lupacciolu [Lu] and C.\ Laurent-Thiebaut [LT], which in our simple case reduces to the fact
that if a continuous function $f$ is $CR$ on
$\Sigma_i$, then for any smoothly bounded domain $D$ in $\Sigma_i$, compactly contained
in $\Sigma_i$ we have
$\int _{D} f\overline\partial \beta = \int_{bD} f\beta $ for
every smooth, $(2,0)$ form on $\C^2$, $i=1, 2$.  The statement in our theorem
is local so we may assume that $E$
is a small perturbation of a piece of a two dimensional plane passing through the center $T$ of $B$
and that the smooth form $\alpha $  has support contained in a small neighbourhood of $T$.
Let $P$ be a small ball centered at $T$ containing the support of $\alpha $ in its interior and,
for small
$\varepsilon >0$, let $P_\varepsilon $
consist of those points of $P$ whose distance from $E$ exceeds $\varepsilon$.
For $i=1,2$ let $P_{\varepsilon, i} = P_\varepsilon\cap\Sigma_i$ and let
$S_{\varepsilon, i} = \{ z\in P\cap\Sigma_i\colon\ dist(z,E)=\varepsilon\}$. The sets
$S_{\varepsilon, i}$ are the only parts of
the boundaries of $P_{\varepsilon_i},\ i=1,2 $ which
interect support of $\alpha$ and as $\varepsilon$ tends to zero, they, as
oriented pieces of the boundaries of $P_{\varepsilon, i},\ i=1,2$, converge
to $E\cap P$ with the opposite orientations.

Now,
$$\eqalign{
\int_M f\overline\partial \alpha &= \lim _{\varepsilon\searrow 0} \int_
{P_\varepsilon}f\overline\partial \alpha \cr
&= \lim _{\varepsilon\searrow 0}
\Bigl[ \int_
{P_{\varepsilon ,1}} f\overline\partial \alpha +\int_
{P_{\varepsilon ,2}} f\overline\partial \alpha \Bigr] \cr
&=\lim _{\varepsilon\searrow 0}
\Bigl[ \int_
{bP_{\varepsilon ,1}} f \alpha +\int_
{bP_{\varepsilon ,2}} f \alpha \Bigr] \cr
&=\lim _{\varepsilon\searrow 0}
\Bigl[ \int_
{S_{\varepsilon, 1}} f \alpha +\int_
{S _{\varepsilon ,2}} f \alpha \Bigr] \cr
&=0\cr}
$$
This completes the proof.
\vskip 2mm
\noindent\bf LEMMA 4.2\ \it Suppose that $U$ is a small open disc whose closure is contained
in $\cS\setminus \R$ and that $G$ is a continuous function on $\pi^{-1}(U)\cap M$ which
is holomorphic on each holomorphic leaf of $\pi ^{-1}(U)\cap N$ and on each holomorphic leaf of
$\pi^{-1}(U)\cap L$. Then $G$ is CR on $\pi^{-1}(U)\cap M$, that is
$$
\int_{\pi^{-1}(U)\cap M} G \overline\partial \beta = 0
$$
for each smooth two zero form $\beta $ on $\C^2$ whose
support meets $\pi^{-1}(U)\cap M$ in a compact set. \rm
\vskip 2mm
The lemma says that if $G$ is CR on $[\pi^{-1}(U)\cap L]\setminus (L\cap N)$ and on
$[\pi^{-1}(U)\cap N]\setminus (L\cap N)$ then $G$ is CR
on $\pi^{-1}(U)\cap M$. This obviously follows from Lemma 4.1 and the
fact that if a continuous function is holomorphic on each holomorphic leaf then it is CR.
\vskip 2mm
\noindent \bf LEMMA 4.3\ \it Let $U\subset\cS\setminus\R$ be an
open disc and let $M_U =\pi^{-1}(U)\cap M$. Suppose that $G$
is a continuous CR function on $M_U$, that is, given a smooth two-zero
form $\omega $ whose support intersects $M_U$ in
a compact set, we have
$$
\int _{M_U}G\overline\partial \omega = 0.
$$
Then the function $z\mapsto \int _{M_z} G(z,w)dw $ is holomorphic on $U$. \rm
\vskip 2mm
\noindent \bf Proof.\ \rm Note first that if $K\subset U$ is a compact set then
$\pi^{-1}(K)$ intersects $M_U$
in a compact set. Let $\alpha $ be a smooth function of one complex
variable $z$ with
compact support contained in $U$. Then
$\beta (z,w) = \alpha (z)  dz\wedge dw $ is a smooth form on $\C^2$
whose support intersects $M_U$ in a compact set so
$$
\int _M G(z,w) {{\partial\alpha (z)}
\over{\partial\overline z}} d\overline z\wedge dz\wedge dw = 0
$$
which, by Fubini, implies that
$$
\int_U\Bigl(\int_{M_z}G(z,w)dw\Bigr){{\partial \alpha}\over{\partial\overline z}}dz = 0.
$$
Since this holds for every smooth function  $\alpha $ with compact support contained in $U$ 
it follows that the function
$z\mapsto \int_{M_z} G(z,w)dw$ is holomorphic on $U$. This completes the proof.
\vskip 2mm
We now show that for each $z\in\cS\setminus \R$ the function $w\mapsto F(z,w)$,
defined on $M_z$, extends holomorphically through $D_z$. Consider the function
$$
H(z,w,W)= {{F(z,w)}\over{w-W}}.
$$
For each $\eta >0$ there is an $R(\eta)<\infty$ such that if $|W|>R(\eta )$, the function
$(z,w)\mapsto H(z,w,W)$ is well defined and continuous
on $P_\eta = \{(z,w)\in M\colon z\in\cS_\eta\}$ where $\cS_\eta = z\in\cS,
\
|\Im z|>\eta\}$, and is holomorphic on each holomorphic leaf in $P_\eta $ so
by Lemma 4.2 it is
CR on $P_\eta $. Lemma 4.3 now implies that for each fixed $W, |W|>R$,
$$
z\mapsto {1\over{2\pi i}}\int_{M_z}{{F(z,w)}\over{w-W}}dw = \Theta (z,W)
$$
is holomorphic on $P_\eta$. Since $z\mapsto \Theta (z,W)$ is continuous on
$\overline{\cS_\eta}$ and since the curves $M_z$ shrink
to $\overline a$ when $z\in\cS_\eta$ approaches
$a\in\bD\setminus\R$ it follows that $\Theta (z,W)=0$ for each $z\in\cS_\eta\cap
\bD\setminus\R$ and thus $\Theta (z,W)\equiv 0\ (z\in\cS_\eta)$. Thus, for each $z\in \cS_\eta $ we
have $\Theta (z,W)=0$ for all $W, |W|>R(\eta)$ which implies that
$$
{1\over{2\pi i}}\int_{M_z} {{F(z,w)}\over{w-W}}dw = 0
$$
for all $W\in \C\setminus\overline{D_z}$ which implies that
the function $w\mapsto F(z,w)$ extends from $M_z$ holomorphically
through $D_z$ for each $z\in\cS_\eta$. Since $\eta>0$
was arbitrary this holds for each $z\in\cS\setminus\R$.

Recall that $(-1+2\tau,0)\times\{\infty\}\subset M$ and that $F$ is continuous on $M$.
Given $T$,\ $-1+2\tau<T<0$,  we will show that $F$ is constant on $\{ T\}\times M_T$.
To do this, we use the reasoning of Tumanov. Fix $T\in (-1+2\tau, 0)$ and observe that
for small $\eta >0$, \ $M_{T+i\eta} $
are simple closed curves bounding $D_{T+i\eta}$ which
depend continuously on $\eta$ and, as domains in $\overline\C$,  continuously tend to the
halfplane $\Im\z <0$ as $\eta $ tends to $0$. Since for each small $\eta >0$ the function $\z\mapsto
F(T+i\eta ,\z )$ extends from $M_{T+i\eta}$ holomorphically through
$D_{T+i\eta}$, the continuity of $F$ implies that $t\mapsto F(T,t)$ has a
bounded holomorphic extension from $\R $ trough  the halfplane $\Im\z <0$.
Repeating the reasoning with $\eta <0$ we
see that $t\mapsto F(T,t)$ has a bounded holomorphic extension from $\R$
through the upper halfplane. Thus, $t\mapsto F(T,t)$ has a bounded holomorphic extension
through $\C$ which, by the
Liouville theorem, must be constant. Thus, for each $T,\ -1+2\tau<T<0$,
the holomorphic extensions of $f$
from all circles in our family which surround $T$, coincide. This implies that $f$ is 
holomorphic in a neighbourhood of
the segment $(-1+2\tau,0)$ and it is easy to see that the analyticity propagates along 
the circles
so it follows that $f$ is holomorphic on $\D$. This
completes the proof of Theorem 1.1.
\vskip 4mm
\bf 5.\ Remarks \rm
\vskip 2mm
A careful examination of the proof of Theorem 1.1 shows that to prove that $f$ is
holomorphic there is no need to assume that $f$ extends holomorphically from
each circle centered at the origin. In fact, the same proof gives
\vskip 2mm
\noindent\bf THEOREM 5.1\ \it Let $f$ be a continuous function on $\DD$ and let $p\in \bD,\ 0<r<1, 0<\rho<1$. Assume that
$f$ extends holomorphically from each circle of radius $R,\ r\leq R\leq 1$, centered at 
the origin, and from
each circle of radius $R,\ \rho\leq R\leq 1$, passing through $p$ and contained in $\DD$ 
If the smallest circles of
these two families are disjoint then $f$ is holomorphic on $\D$.\rm
\vskip 2mm
Let $p_1, p_2\in b\D,\ p_1\not=p_2$. In a way similar to the way above
we prove that a continuous function on $\DD$
which extends holomorphically from each circle contained in $\DD$ and
passing through $p_1$ and which extends holomorphically from each circle contained in $\DD$
and passing through $p_2$ then $f$ is holomorphic on $\D$. In fact, again, fewer circles suffice:
\vskip 2mm
\noindent\bf THEOREM 5.2\ \it Let $p_1, p_2\in\bD,\ p_1\not=p_2$. Let $0<r_j<1,\ j=1,2$,  and assume that that $f$
is a continuous function on $\DD$ which, for each $j=1,2$, extends holomorphically from each circle of radius
$\rho,\ r_j\leq\rho\leq 1$, contained in $\DD$ and passing
through $p_j$. If the smallest circles of these two families are disjoint then $f$ is holomorphic on $\D$. \rm
\vskip 2mm
\noindent Note that, after applying an automorphism of $\D$ one can, with no loss of
generality assume that $p_1=-1,\ p_2=1$.
The domains $D_z$ now are bounded by two circular arcs.  Note that the example in
Section 1 shows that in both theorems
the condition that the smallest circles
of the families be disjoint cannot be omitted.
\vskip 5mm
\noindent\bf ACKNOWLEDGEMENT \ \rm The author is indebted to E.\ L.\ Stout for the proof of Lemma 4.1.

This work was supported
in part by the Ministry of Higher Education, Science and Technology of Slovenia
through the research program Analysis and Geometry, Contract No.\ P1-0291(B)
\vskip 10mm
\centerline{\bf References}
\vskip 2mm

\noindent [A]\ M.\ L.\ Agranovsky:  Propagation of boundary CR foliations and Morera type
theorems for manifolds with attached analytic discs.

\noindent Adv. Math. 211 (2007) 284--326.
\vskip 2mm
\noindent [AG]\ M.\ L.\ Agranovsky and J.\ Globevnik: Analyticity
on circles for rational and real-analytic functions of two real variables.

\noindent J.\ d'Analyse Math.\ 91 (2003) 31-65
\vskip 2mm
\noindent [G1] \ J.\ Globevnik: Holomorphic extensions from open families of circles.

\noindent Trans.\ Amer.\ Math.\ Soc.\ 355 (2003) 1921-1931
\vskip 2mm
\noindent [G2]\ J.\ Globevnik: Analyticity on families of circles.

\noindent Israel J.\ Math. 142 (2004) 29-45
\vskip 2mm
\noindent [G3] \ J.\ Globevnik: Analyticity on translates of a Jordan curve.

\noindent Trans.\ Amer.\ Math.\ Soc. 359 (2007) 5555-5565
\vskip 2mm
\noindent [L]\ H.\ Lewy: On the local character of the
solutions of an atypical linear differential equation in three
variables and a related theorem for regular functions of two complex variables.

\noindent Ann. of Math.\ 64 (1956) 514-522
\vskip 2mm
\noindent [Lu]\ G.\ Lupacciolu:\ A theorem on holomorphic extensions of CR-functions.

\noindent Pacif.\ J.\ Math.\ 124 (1986) 177-191
\vskip 2mm
\noindent [LT]\ C.\ Laurent-Thiebaut:\ Sur l'extension des fonctions CR dans
une variete de Stein.

\noindent Ann.\ Mat.\ Pura Appl.\ 150 (1988) 141--151
\vskip 2mm
\noindent [R]\ H.\ Rossi: A generalization of a theorem of Hans Lewy.

\noindent Proc.\ Amer.\ Math.\ Soc.\ 19 (1968) 436-440
\vskip 2mm
\noindent [T1]\ A.\ Tumanov: A Morera type theorem in the strip.

\noindent Math.\ Res.\ Lett.\ 11 (2004) 23-29
\vskip 2mm
\noindent [T2]\ A.\ Tumanov: Testing analyticity on circles.

\noindent Amer.\ J.\ Math.\ 129 (2007) 785-790

\vskip 8mm
\noindent Institute of Mathematics, Physics and Mechanics

\noindent University of Ljubljana, Ljubljana, Slovenia

\noindent josip.globevnik@fmf.uni-lj.si

\bye